\documentclass{amsart}
\usepackage{amssymb,amsmath,amsthm}
\usepackage[utf8]{inputenc}
\usepackage[italian,english]{babel}
\usepackage{mathrsfs}
\usepackage[textmath,displaymath,floats,graphics]{preview}
\usepackage[all]{xy}
\usepackage{pdfsync}
\usepackage{mathabx}
\usepackage{tikz}
\usetikzlibrary{matrix,arrows}

\theoremstyle{definition}
\newtheorem{definition}{Definition}
\theoremstyle{remark}
\newtheorem{example}{Example}
\newtheorem{remark}{Remark}
\theoremstyle{plain}
\newtheorem{theorem}{Theorem}
\newtheorem{proposition}{Proposition}
\newtheorem{lemma}{Lemma}
\newtheorem{corollary}{Corollary}

\newenvironment{sistema}
{\left\lbrace\begin{array}{@{}l@{}}}
{\end{array}\right.}

\begin{document}

\markboth{Paola Bonacini}{On the lifting problem in $\mathbb P^4$ in characteristic $p$}

\title{On the lifting problem in $\mathbb P^4$ in characteristic $p$}
\author{Paola Bonacini}
\address{Dipartimento di Matematica e Informatica, Università di Catania\\
Viale A. Doria 6, 95124, Catania, Italy\\ 
bonacini@dmi.unict.it} 

\maketitle

\begin{abstract}
Given $\mathbb P^4_k$, with $k$ algebraically closed field of
characteristic $p>0$, and $X\subset \mathbb P^4_k$ integral surface of
degree $d$, let $Y=X\cap H$ be the general hyperplane section of $X$. We suppose
that $h^0\mathscr I_Y(s)\ne 0$ and $h^0\mathscr I_X(s)=0$ for some
$s>0$. This determines a nonzero element
$\alpha\in H^1\mathscr I_X(s)$ such that $\alpha\cdot H=0$ in
$H^1\mathscr I_X(s)$. We find different upper bounds of $d$ in terms of
$s$, $p$ and the order of $\alpha$ and we show that these bounds are
sharp. In particular, we see that $d\le s^2$ for $p<s$ and $d\le
s^2-s+2$ for $p\ge s$.
\end{abstract}

%\keywords{surface, sporadic zero, plane curves}

%\ccode{Mathematics Subject Classification 2000: 14H50, 14M07, 14J99}

\section{Introduction}

Let $X\subset \mathbb P^4_k$, with $k$ algebraically closed field of
positive characteristic $p$, be
an integral surface of degree $d$. Let $Y=X\cap H$ be the general
hyperplane section of $X$ and consider a surface of degree $s$
containing $Y$. In this paper we study the problem of lifting a surface of $H$ of
degree $s$ containing $Y$ to a hypersurface in $\mathbb P^4$ of degree $s$ containing
$X$.  In particular, we suppose that $h^0\mathscr I_Y(s)\ne 0$
and $h^0\mathscr I_X(s)=0$ for some $s>0$. 

In the case that $\operatorname{char}k=0$ the problem
has been studied and solved by Mezzetti and Raspanti in \cite{MeR} and in \cite{Me1}, showing that
$d\le s^2-s+2$ and that this bound is sharp, and in \cite{Me} Mezzetti
classifies the border case $d=s^2-s+2$. Other results concerning the
lifting problem have been obtained in characteristic $0$ for curves in
$\mathbb P^3$ (see \cite[Corollary p.~147]{L}, \cite{GP} and
\cite[Corollario~2]{S1}) and for integral varieties of codimension $2$
in $\mathbb P^n$ (see, for example, \cite{Me1} for $n=5$, \cite{ROG2}
for $n=6$ and \cite{ROG5} for the general case). In the case that
$\operatorname{char}p>0$ the lifting problem has been studied for
curves in $\mathbb P^3$ in \cite{BON}.

In this paper, the starting point is that the non lifting section of $H^0\mathscr I_Y(s)$ determines a nonzero element
$\alpha\in H^1\mathscr I_X(s)$ such that $\alpha\cdot H=0$ in
$H^1\mathscr I_X(s)$. The order of $\alpha$ is the maximum integer
$m\in \mathbb N$ such that $\alpha=\beta\cdot H^m$ for some $\beta\in
H^1\mathscr I_X(s-m-1)$. For $p<s$ we need to relate $s$ and $m$. In
particular, in Theorem~\ref{T} we suppose that $p<s$ and
we show that: 
\begin{enumerate}
\item $d\le s^2-s+2+p^n$, if $s\ge 2m+3$, with $p^n\le m+1$ and
  $p^{n+1}>m+1$;
\item $d\le s^2$ if
$s\le 2m+2$.
\end{enumerate}
As a consequence, we see that for $p<s$ it must be $d\le s^2$. For
$p\ge s$ in Theorem~\ref{T:2} we show that $d\le s^2-s+2$, which is
the same bound as in the characteristic $0$ case. In Example~\ref{ex}
we see that the bounds given in Theorem~\ref{T} and in Theorem~\ref{T:2} are sharp.

\section{Hilbert function of points in $\mathbb P^2$}

Let us denote by $X$ a zero-dimensional scheme in $\mathbb P^2_k$,
where $k$ is an algebraically closed field of any characteristic. Let $H_X\colon \mathbb
N\rightarrow \mathbb N$ be the Hilbert function of $X$ and let us
consider the first difference of $H_X$:
\[\Delta H(X,i)=H(X,i)-H(X,i-1).\]
It is known \cite{GMP} that there exist $a_1\le a_2\le t$ such that:
\[
\Delta H(X,i)=
 \begin{cases}
    i+1                   & \text{for $i=0,\dots,a_1-1$} \\
    a_1                   & \text{for $i=a_1,\dots,a_2-1$} \\
    <a_1                  & \text{for $i=a_2$} \\
    \text{non increasing} & \text{for $i=a_2+1,\dots,t$} \\
    0                     & \text{for $i>t$.}
 \end{cases}
\]
\begin{definition}
We say that $X$ has \emph{the Hilbert function of decreasing type} if for $a_2\le i<j<t$ we have $\Delta H(X,i)>\Delta H(X,j)$.
\end{definition}

The following is a result well known in characteristic $0$ (see
\cite{HAR} and \cite[Corollary~2]{MaRa3}) and proved in any
characteristic in \cite[Corollary~4.3]{BON2}.
\begin{theorem} \label{T:3}
Let $C\subset \mathbb P^3$ be an integral curve and let $X$ be its
general plane section. Then $H_X$ is of decreasing type.
\end{theorem}

The following result will be useful in the proof of the main results of the paper.
\begin{proposition} \label{P:2}
  Let $X\subset \mathbb P^2$ be a $0$-dimensional scheme whose Hilbert
  function is of decreasing type. Let us suppose that $h^0\mathscr
  I_X(s-1)=0$ for some $s>0$ and that one of the following conditions
  holds:
  \begin{enumerate}
  \item $h^0\mathscr I_X(s)\ge 3$;
  \item $h^0\mathscr I_X(s)=2$ and there exists $i\in \mathbb N$ such that $\Delta H_X(s+i)\le s-i-2$.
  \end{enumerate}
Then $\deg X\le s^2-s+i+1$.
\end{proposition}
\begin{proof}
  The proof is a straightforward computation and follows by the fact
  that the Hilbert function of $X$ is of decreasing type. 

If  $h^0\mathscr I_X(s)\ge 3$, then we see that:
\[
\deg X\le \frac{s(s+1)}{2}+\frac{(s-2)(s-1)}{2}=s^2-s+1<s^2-s+i+1.
\]
Let $h^0\mathscr I_X(s)=2$ and let us suppose that $i=\min\{k\in \mathbb
N\mid \Delta H_X(s+k)\le s-k-2\}$. Since $H_X$ is of decreasing type,
$\Delta H_X(s+k)=s-k-1$ for $k\le i-1$ and $\Delta H_X(s+k)\le s-k-2$
for $k\ge i$. Then: 
\[
\deg X\le
\frac{s(s+1)}{2}+\sum_{k=0}^{i-1}(s-k-1)+\sum_{k=i}^{s-3}(s-k-2)=s^2-s+i+1.
\]
\end{proof}

\section{Frobenius morphism and incidence varieties}

In this section we show some results about incidence varieties and
Frobenius morphism. First let us recall the definition of absolute and relative Frobenius morphism:

\begin{definition}
The absolute Frobenius morphism of a scheme $X$ of characteristic $p>0$ is $F_X\colon X\rightarrow X$, where $F_X$ is the identity as a map of topological spaces and on each $U$ open set $F^{\#}_{X}\colon\mathscr O_X(U)\rightarrow \mathscr O_X(U)$ is given by $f\mapsto f^p$ for each $f\in \mathscr O_X(U)$. Given $X\rightarrow S$ for some scheme $S$ and $X^{p/S}=X\times_{S,\, F_S}S$, the absolute Frobenius morphisms on $X$ and $S$ induce a morphism $F_{X/S}\colon X\rightarrow X^{p/S}$, called the Frobenius morphism of $X$ relative to $S$.
\end{definition}

Given $\mathbb P^n$ for some $n\in \mathbb N$, let us consider the
bi-projective space $\widecheck{\mathbb P}^n\times\mathbb P^n$ and let
$r\in \mathbb N$ be a non negative integer. Let $k[\underline t]$ and
$k[\underline x]$ be the coordinate rings of $\widecheck{\mathbb P}^n$
and $\mathbb P^n$, respectively. Let $M_r\subset \widecheck{\mathbb P}^n\times\mathbb P^n$ be the hypersurface of equation:
\[h_r:=\sum_{i=0}^nt_i{x_i}^{p^r}=0.\]
Note that in the case $r=0$ $M_r$ is the usual incidence variety $M$ of equation $\sum t_ix_i=0$. If $r\ge 1$, $M_r$ is determined by the following fibred product:
\begin{equation}
  \label{eq:2}
  \begin{tikzpicture}[baseline=(current  bounding  box.center), descr/.style={fill=white,inner sep=2.5pt}]
\matrix(m)[matrix of math nodes, row sep=2.6em, column sep=2.8em, text
height=1.5ex, text depth=0.25ex] {M&&\\& M_r&M\\ & \mathbb P^n &
  \mathbb P^n\\}; 
\path[dashed, ->] (m-1-1) edge node[descr] {$\scriptstyle{F_{M_r}}$} (m-2-2);
\path[->] (m-1-1) edge [bend right=20] node[auto, left] {$\scriptstyle
  p$} (m-3-2);
\path[->] (m-1-1) edge [bend left=20] node[auto, right] {$\scriptstyle{(F_M)^r}$} (m-2-3);
\path[->] (m-2-2) edge node[auto] {$\scriptstyle \pi$} (m-2-3); 
\path[->] (m-2-2) edge node[left] {$\scriptstyle{p_{M_r}}$} (m-3-2); 
\path[->] (m-2-3) edge node[auto] {$\scriptstyle{p}$} (m-3-3); 
\path[->] (m-3-2) edge node[below] {$\scriptstyle{F^r}$} (m-3-3); 
\end{tikzpicture}
\end{equation}
where $F\colon \mathbb P^n\rightarrow\mathbb P^n$ is the absolute
Frobenius. 

\begin{remark}
  \label{R:1}
$M=\mathbb P(\mathscr T_{\mathbb P^n}(-1))$, so that by
\cite[Lemma~1.5]{EIN} $M_r=\mathbb P({F^r}^{\star}(\mathscr T_{\mathbb
  P^n}(-1)))$ and by \cite[Ch.II, ex.~7.9]{HART} we see that $\operatorname{Pic}(M_r)=\mathbb Z\times \mathbb Z$
for any $r\ge 0$. Moreover, since we have:
\[
0\rightarrow \mathscr O_{\widecheck{\mathbb P}^n\times\mathbb
  P^n}(-1,-p^r)\rightarrow \mathscr O_{\widecheck{\mathbb P}^n\times\mathbb
  P^n}\rightarrow \mathscr O_{M_r}\rightarrow 0
\]
and $H^1 \mathscr O_{\widecheck{\mathbb P}^n\times\mathbb
  P^n}(m,n)=0$ for any $m,n\in \mathbb Z$ by the K\"unneth formula
\cite[Ch.VI, Corollary~8.13]{M}, then any hypersurface $V\subset M_r$
is the complete intersection given by $g=h_r=0$ for some
bi-homogeneous $g\in k[\underline t,\underline x]$.
\end{remark}

Let $\eta\in \widecheck{\mathbb P}^n$ be the generic point and and
consider $g_{M_r}\colon M_r\rightarrow \widecheck{\mathbb P}^n$. Then
$g_{M_r}^{-1}(\eta)$ is isomorphic to the hypersurface $H_r$ of $\mathbb P^n$ of
degree $p^r$ such that, over the algebraic closure
$\overline{k(\eta)}$ of $k(\eta)$,
$({H_r})_{\text{red}}$ is the generic hyperplane $H$ of $\mathbb P^n$.

\begin{proposition}   \label{P:6}
 $\Omega_{M_r/\mathbb P^n}|_H\cong {F^r}^{\star}\mathscr T_H(-p^r)$.
\end{proposition}
\begin{proof}
  The sheaf $\mathscr
  E={F^r}^{\star}(\mathscr T_{\mathbb P^n}(-1))={F^r}^{\star}\mathscr
  T_{\mathbb P^n}(-p^r)$ is determined by the exact sequence
  $0\rightarrow \mathscr O_{\mathbb P^n}(-p^r)\rightarrow {\mathscr
    O_{\mathbb P^n}}^{\oplus n+1}\rightarrow \mathscr E\rightarrow 0$
  and, since $M_r=\mathbb
  P({F^r}^{\star}(\mathscr T_{\mathbb P^n}(-1)))$, by \cite[Ch.III, Ex.~8.4(b)]{HART} we have also $0\rightarrow \Omega_{M_r/\mathbb P^n}\rightarrow {p_{M_r}}^{\star}({F^r}^{\star}\mathscr T_{\mathbb P^n}(-p^r))\otimes_{\mathscr O_{M_r}}\mathscr O_{M_r}(-1,0)\rightarrow \mathscr O_{M_r}\rightarrow 0$. When we restrict to $H$, by the fact that the sequence locally splits it follows that the following sequence is exact:
\[0\rightarrow \Omega_{M_r/\mathbb P^n}|_H\rightarrow {p_{M_r}}^{\star}({F^r}^{\star}\mathscr T_{\mathbb P^n}(-p^r)) \otimes_{\mathscr O_{M_r}}\mathscr O_{M_r}(-1,0)|_H\rightarrow \mathscr O_H\rightarrow 0\]
\begin{equation}
  \label{eq:1}
\Rightarrow  0\rightarrow \Omega_{M_r/\mathbb P^n}|_H\rightarrow {F^r}^{\star}\mathscr T_{\mathbb P^n}(-p^r)|_H\rightarrow \mathscr O_H\rightarrow 0.
\end{equation}
Since $\mathscr T_{\mathbb P^n}(-1)|_H\cong \mathscr T_H(-1)\oplus
\mathscr O_H$, then ${F^r}^{\star}\mathscr T_{\mathbb
  P^n}(-p^r)|_H\cong {F^r}^{\star}\mathscr T_H(-p^r)\oplus\mathscr
O_H$. By the fact that ${F^r}^{\star}\mathscr T_H$ is stable (see
\cite[Ch.II, Theorem~1.3.2]{OSS} and \cite[Theorem~2.1]{MeRa}) and that $\mu({F^r}^{\star}\mathscr T_H(-p^r))>0$, it follows that $\operatorname{Hom}({F^r}^{\star}\mathscr T_H(-p^r),\mathscr O_H)=0$ and so by \eqref{eq:1} also that $\Omega_{M_r/\mathbb P^n}|_H\cong {F^r}^{\star}\mathscr T_H(-p^r)$.
\end{proof}

Now we prove some results about the projection from a hypersurface in
$M_r$ to $\mathbb P^n$. 

\begin{theorem}  \label{T:1}
Let $V\subset \widecheck{\mathbb P}^n\times \mathbb P^n$ be an integral hypersurface in $M$ such that the projection $\pi\colon V\rightarrow \mathbb P^n$ is dominant and not generically smooth. Then there exist $r\ge 1$ and $V_r\subset M_r$ integral hypersurface such that $\pi$ can be factored in the following way:
\begin{center}
\begin{tikzpicture}
\matrix(m)[matrix of math nodes, row sep=2.6em, column sep=2.8em, text
height=1.5ex, text depth=0.25ex] {V& &\mathbb P^n\\ & V_r &\\}; 
\path[->] (m-1-1) edge node[auto] {$\scriptstyle \pi$} (m-1-3);
\path[->] (m-1-1) edge node[auto, left] {$\scriptstyle F_r$} (m-2-2) ;
\path[->] (m-2-2) edge node[auto, right] {$\scriptstyle \pi_r$} (m-1-3); 
\end{tikzpicture}
\end{center}
where the projection $\pi_r$ is dominant and generically smooth and $F_r$ is induced by the commutative diagram:
\begin{center}
\begin{tikzpicture}
\matrix(m)[matrix of math nodes, row sep=2.6em, column sep=2.8em, text
height=1.5ex, text depth=0.25ex] {V&V_r\\ M&M_r\\}; 
\path[->] (m-1-1) edge node[auto] {$\scriptstyle F_r$} (m-1-2);
\path[right hook->] (m-1-1) edge node[auto, left] {$\scriptstyle j$} (m-2-1) ;
\path[right hook->] (m-1-2) edge node[auto] {$\scriptstyle i$} (m-2-2); 
\path[->] (m-2-1) edge node[below] {$\scriptstyle{F_{M_r}}$} (m-2-2); 
\end{tikzpicture}
\end{center}
\end{theorem}
\begin{proof}
 The proof works as in \cite[Theorem~3.3]{BON}. Indeed, the proofs of Lemma 3.1 and
 Proposition 3.2 in \cite{BON} works also for $\widecheck{\mathbb P}^n\times\mathbb
  P^n$.
\end{proof}

\begin{proposition}
  \label{P:4}
  Let $V_r\subset M_r$ be an integral hypersurface given by:
\[
\begin{sistema}
  q(\underline t,\underline x)=0\\[2ex]
\displaystyle\sum_{i=0}^n {t_i}{x_i}^{p^r}=0
\end{sistema}
\] 
such that the
  projection $\pi_r\colon V_r\rightarrow \mathbb P^n$ is generically
  smooth. Then $\pi_r$ is not smooth exactly on the following closed subset
  of $V_r$:
\[
V_r\cap V\left({x_i}^{p^r}\frac{\partial q}{\partial
    t_i}-{x_j}^{p^r}\frac{\partial q}{\partial t_i} \mid i,j=0,\dots,n\right).
\]
\end{proposition}
\begin{proof}
  Let $P_0=(\underline a,\underline b)\in V_r$ be such that $V_r$ is
  not smooth in $P_0$. Then there exists $\lambda\in k$ such that:
  \begin{equation}
    \label{eq:3}
    \dfrac{\partial q}{\partial t_i}(P_0)=\lambda {b_i}^{p^r}
  \end{equation}
for any $i=0,\dots,n$. 

If $P_0$ is a regular point, then the projective tangent space $T_{V_r,P_0}$ at $P_0\in V_r$ is given by the equations:
\[
\sum_{i=0}^n \frac{\partial q}{\partial x_i}(P_0)x_i+\sum_{i=0}^3\frac{\partial q}{\partial t_i}(P_0)t_i=\sum_{i=0}^n(a_ix_i+b_it_i)=0\]
if $r=0$ and by the equations:
\[
\sum_{i=0}^n \frac{\partial q}{\partial x_i}(P_0)x_i+\sum_{i=0}^n\frac{\partial q}{\partial t_i}(P_0)t_i=\sum_{i=0}^n{b_i}^{p^r}t_i=0
\]
if $r\ge 1$. In both cases the projection on $T_{\mathbb P^n\!,{\pi}(P_0)}$ is not surjective if and only if there exists $\lambda\in k$ such that:
\[\dfrac{\partial q}{\partial t_i}(P_0)=\lambda {b_i}^{p^r}\quad
\forall\, i=0,\dots,n.\]
This together with \eqref{eq:3} proves the statement. 
\end{proof}

\section{Lifting problem}
Let $X\subset \mathbb P^4$ be a scheme and, following the previous
notation, consider the projections $p_{M_r}\colon M_r\rightarrow \mathbb P^4$ and $g_{M_r}\colon M_r\rightarrow \widecheck{\mathbb P}^4$. Let $T_r={p_{M_r}}^{-1}(X)$ and:
\[\mathscr{I}_r(m,n)={g_{M_r}}^{\star}\left(\mathscr O_{\check{\mathbb P^4}}(m)\right)\otimes_{\mathscr O_{M_r}}{p_{M_r}}^{\star}\left(\mathscr I_X(n)\right)\]
for every $m$, $n\in \mathbb Z$.
\begin{proposition}  \label{P:1}
If $\mathscr I_r=\mathscr I_r(0,0)$ and $\mathscr I_{T_r}$ is the ideal sheaf of $T_r$ in $M_r$, then $\mathscr I_r=\mathscr I_{T_r}$.
\end{proposition}
\begin{proof}
  The proof works as in \cite[Proposition 4.1]{BON}.
\end{proof}

Let $X\subset \mathbb P^4$ be an integral surface of degree $d$. Let $Y=X\cap H$ be the generic
hyperplane section of $X$ and let $Z=Y\cap K$ be the generic plane section
of $Y$. Let $\mathscr I_X$ be the ideal sheaf of $X$ in
$\mathbb P^4$,  $\mathscr I_Y$ the ideal sheaf of $Y$ in
$H\cong \mathbb P^3$ and $\mathscr I_Z$ the ideal sheaf of $Z$ in
$K\cong \mathbb P^2$. Let us consider for any $s\in \mathbb N$ the
following maps:
\[
\pi_s\colon H^0\mathscr I_X(s)\rightarrow
H^0\mathscr I_Y(s) \mbox{\quad and \quad}\phi_s\colon H^1\mathscr
I_X(s-1)\rightarrow H^1\mathscr I_X(s)
\]

A \emph{sporadic zero} of degree $s$ is an element $\alpha\in \operatorname{coker} (\pi_s)=\operatorname{ker} (\phi_s)$.

\begin{definition}
  \label{D:1}
  The \emph{order} of a sporadic zero $\alpha$ is the maximum integer $m$ such that
  $\alpha=\beta\cdot H^m$, for some $\beta \in H^1\mathscr
  I_X(s-m-1)$.
\end{definition}

\begin{proposition}
  \label{P:5}
  Let $\alpha$ be a sporadic zero of degree $s$ and let $h^0\mathscr
  I_X(s)=0$. Then one of the following conditions holds:
  \begin{enumerate}
  \item $\deg X\le s^2-s+1$;
    \item $h^0\mathscr I_Y(s)=1$ and $h^0\mathscr I_Z(s)=2$.
  \end{enumerate}
\end{proposition}
\begin{proof}
  Let $q=\min\{i\mid h^0\mathscr
I_Y(i)\ne 0\}$. So $q\le s$ and by hypothesis there is an integral
surface of degree $q$ containing $Y$ that does not lift to an integral
surface of degree $q$ containing $X$. In particular we have a sporadic
zero of degree $q$ for $X$ and by \cite[Theorem 2.1]{ROG4} we get a sporadic zero for $Y$ of
degree $s'\le q$. By \cite[Theorem 4.1]{BON2} this
means that there is an integral curve of degree $s'$ in $K$ containing $Z$ that does
not lift to a surface in $H$ of degree $s'$ containing $Y$. However, by
restricting the integral surface of degree $q$ containing $Y$ to $K$
we see that $d\le qs'$. 

If $s'<s$, then we see that
$\deg X=\deg Z\le s^2-s$. 

So we can suppose that $q=s'=s$, which implies that $h^0\mathscr
I_Z(s)\ge 1+h^0\mathscr I_Y(s)\ge 2$. If
$h^0\mathscr I_Z(s)\ge 3$, then by Theorem \ref{T:3} and by Proposition \ref{P:2} we get $\deg
X=\deg Z\le s^2-s+1$. So we can suppose that $h^0\mathscr I_Z(s)=2$,
which implies also that $h^0\mathscr I_Y(s)=1$.
\end{proof}

The following result, together with Proposition \ref{P:5}, provides us with
the tools for the proof of the
main results of this paper.

\begin{lemma}   \label{L:1}
  Let $\alpha$ be a sporadic zero of degree $s$ and order $m$. Suppose
  that $\alpha$ determines a non-liftable integral surface $R$ in $H$ of degree $s$
  containing $Y$ and that $I_R=(f)$ for some $f\in H^0\mathscr
  O_H(s)$. Then for some $r\in \mathbb N$ such that $p^r\le m+1$ there
  exist:
  \begin{enumerate}
  \item $f_i\in H^0\mathscr O_{H}(s)$ for $i=0,\dots,4$ such that the subscheme of $H$ associate to the ideal
    $(f,{x_i}^{p^r}f_j-{x_j}^{p^r}f_i|_H,i,j=0,\dots,4)$ is a
    $1$-dimensional scheme
    $E$ (which can have isolated or embedded $0$-dimensional subschemes) such that $Y\subset E \subset R$;
  \item a reflexive sheaf $\mathscr N$ of rank $3$ such that we have
    the exact sequence:
    \begin{equation}
      \label{eq:11}
      0\rightarrow \mathscr N\rightarrow {F^r}^{\star}\Omega_{H}(p^r)\rightarrow \mathscr I_{E|R}(s)\rightarrow 0,
    \end{equation}
being $\mathscr I_{E|R}\subset \mathscr O_{R}$ the ideal sheaf of $E$.
  \end{enumerate}
\end{lemma}
\begin{proof}
Let $r\ge 0$ and let $\mathscr I_r={p_{M_r}}^{\star} \mathscr
I_X$. Given the generic point $\eta\in
\widecheck{\mathbb P}^4$ and $g_{M_r}\colon M_r\rightarrow \widecheck{\mathbb
P}^4$, we have seen that $g_{M_r}^{-1}(\eta)$ is isomorphic to the hypersurface
$H_r$ of $\mathbb P^4$ of
degree $p^r$ such that, over $\overline{k(\eta)}$, $({H_r})_{\text{red}}=H$.

By proceeding as in \cite[Theorem 1.2, Step 1 and Step 2]{BON} and by Theorem \ref{T:1} we see
that there exist $r\ge 0$ and  $V_r\subset M_r$ hypersurface given by:
\[
\begin{sistema}
   q(\underline t,\underline x)=0\\[2ex]
\displaystyle\sum_{i=0}^4 {t_i}{x_i}^{p^r}=0
\end{sistema}
\]
such that the projection $p_{V_r}\colon V_r\rightarrow \mathbb P^4$ is
generically smooth and, given $g_{V_r}\colon V_r\rightarrow \widecheck {\mathbb P}^4$,
$g_{V_r}^{-1}(\eta)$ is the complete intersection of $H_r$
with a hypersurface of $\mathbb P^4$ of degree $s$ and it is such that
$g_{V_r}^{-1}(\eta)_{\text{red}}\cong R$ over
$\overline{k(\eta)}$. This means that $m\ge p^r-1$.

Let $U\subset V_r$ be the subscheme where $p_{V_r}$ is not
smooth. Then by Proposition \ref{P:4} we see that:
\begin{equation}
  \label{eq:9}
  U=V_r\cap V\left({x_i}^{p^r}\frac{\partial q}{\partial
    t_i}-{x_j}^{p^r}\frac{\partial q}{\partial t_i} \mid i,j=0,\dots,4\right).
\end{equation}
By proceeding as in \cite[Theorem 1.2, Step 3]{BON} we see that
$U\supseteq T_r$, $\dim U=5$ and we have for some $b>0$:
\begin{equation}  \label{E:2}
0\rightarrow {\Omega_{V_r/\mathbb P^4}}^{\vee}\rightarrow {\Omega_{M_r/\mathbb P^4}}^{\vee}\otimes_{\mathscr O_{M_r}} \mathscr O_{V_r}\rightarrow \mathscr I_{U|V_r}(b,s)\rightarrow 0,
\end{equation}
with $\mathscr I_{U|V_r}\subset \mathscr O_{V_r}$ ideal sheaf of
$U$. 

Restricting \eqref{E:2} to $H$, by Proposition \ref{P:6} we get a surjective map
${F^r}^{\star}\Omega_{H}(p^r)\otimes_{\mathscr O_{H}} O_{R}\rightarrow \mathscr I_{E|R}(s)$, with $\mathscr
I_{E|R}\subset \mathscr O_{R}$ ideal sheaf of the $1$-dimensional
scheme $E=U\cap g_{M_r}^{-1}(\eta)_{\text{red}}$. Note that $E\supseteq
T_r\cap g_{M_r}^{-1}(\eta)_{\text{red}}\cong Y$. The kernel of the map ${F^r}^{\star}\Omega_{H}(p^r)\rightarrow
\mathscr I_{E|R}(s)$ is the  sheaf $\mathscr N$  that
determines the exact sequence \eqref{eq:11} and it is torsion free and
normal and so it is reflexive. Moreover, by \eqref{eq:9} we get
\[
E=V\left(q|_H,{x_i}^{p^r}\frac{\partial q}{\partial
    t_i}-{x_j}^{p^r}\frac{\partial q}{\partial t_i}|_H \mid i,j=0,\dots,4\right),
\]
where $q|_H=f$, and so the statement is proved by taking $f_i=\frac{\partial q}{\partial
    t_i}|_H$ for any $i=0,\dots,4$. 
  \end{proof}

Now we can prove the first main result of the paper.

\begin{theorem}  \label{T}
Let $\alpha$ be a sporadic zero of degree $s$ and order $m$ and let $p<s$. Suppose that
$h^0\mathscr I_X(s)=0$. Then:
\begin{enumerate}
\item if $s\ge 2m+3$, we have $d\le s^2-s+p^n+1$, with $p^n\le m+1$
  and $p^{n+1}>m+1$;
\item if $s\le 2m+2$, we have $d\le s^2$. 
\end{enumerate}
\end{theorem}
\begin{proof}
By Proposition \ref{P:5} we can suppose that $h^0\mathscr I_Y(s)=1$
and $h^0\mathscr I_Z(s)=2$. In particular, if $s\le 2m+2$, we get the
conclusion. So we suppose that $s\ge 2m+3$ and we also see that the
surface $R$ of degree $s$
  containing $Y$ that can not be lifted to a hypersurface of degree
  $s$ containing $X$ is integral. Let $I_R=(f)$ in $H$ be the ideal of $R$.

By Lemma \ref{L:1} we see that there exist
$r\in \mathbb N$ with $p^r\le m+1$ and $f_i\in H^0\mathscr O_{H}(s)$
for $i=0,\dots,4$ such that the subscheme of $H$ associate to the
ideal $(f,{x_i}^{p^r}f_j-{x_j}^{p^r}f_i|_H,i,j=0,\dots,4)$ is a $1$-dimensional scheme $E$, which can have isolated or embedded $0$-dimensional schemes, such that $Y\subset E \subset R$. Moreover, there exists a reflexive sheaf $\mathscr N$ of rank $3$ such that we have the exact sequence:
    \begin{equation}
      \label{eq:12}
      0\rightarrow \mathscr N\rightarrow {F^r}^{\star}\Omega_{H}(p^r)\rightarrow \mathscr I_{E|R}(s)\rightarrow 0,
    \end{equation}
being $\mathscr I_{E|R}\subset \mathscr O_{R}$ the ideal sheaf of
$E$. We want to prove that $d\le s^2-s+1+p^r$. 

Note that $c_1(\mathscr N)=-p^r-s$ and
\begin{equation}
  \label{eq:8}
  c_2(\mathscr N)=s^2+p^rs+p^{2r}-\deg E.
\end{equation}
So we see that if $\mathscr N$ is semistable, by the Bogomolov
inequality for semistable reflexive sheaves and by the fact that $\deg
E\ge \deg Y=\deg X$ we get the statement. So we can suppose that
$\mathscr N$ is unstable. Moreover by Theorem \ref{T:3}and by Proposition \ref{P:2} we can
suppose that $\Delta H_Z(s+i)=s-i-1$ for any $i\le p^r$. Given $g\in
H^0\mathscr O_K(s)$ such that $f|_K$ and $g$ are generators of $I_Z$ in degree $s$, by
\cite[Proposition1.4]{MaRa2} we see that $f|_K$ and $g$ are the only generators of
$I_Z$ in degree $\le s+p^r$. By this assumptions we will get a contradiction.

Restricting \eqref{eq:12} to $K$ we get:
  \begin{equation}
    \label{eq:5}
    0\rightarrow \mathscr N|_K\rightarrow {F^r}^{\star}\Omega_K(p^r)\oplus
\mathscr O_K\rightarrow \mathscr I_{E\cap K|R\cap K}(s)\rightarrow 0.
  \end{equation}
Since $N$ is unstable of rank 3,
${F^r}^{\star}\Omega_H(p^r)$ is stable and $c_1({F^r}^{\star}\Omega_H(p^r))=-p^r<0$, the maximal
destabilizing subsheaf $\mathscr F$ of $\mathscr N$ has rank at most 2 and 
$c_1(\mathscr F)<0$. By \cite[Theorem 3.1]{MARU} we see that $F|_K$ is
still semistable and so it must be $h^0\mathscr N|_K=0$. By
\eqref{eq:5} we see that $h^0\mathscr I_{E\cap K|R\cap K}(s)\ge 1$,
which implies that $h^0\mathscr I_{E\cap K}(s)\ge 2$ and, since $E\cap
K\supseteq Z$ and $h^0\mathscr I_Z(s)=2$, we get that $H^0\mathscr
I_{E\cap K}(s)=2$. Since $R\cap K$ is integral of degree $s$ and
$R\cap K\supset E\cap K$, we see that $\deg(E\cap K)\le s^2$. 

Recall that for any $i,j=0,\dots,4$:
\[
{x_i}^{p^r}f_j-{x_j}^{p^r}f_i|_H\in H^0\mathscr
I_{E}(s+p^r)\Rightarrow {x_i}^{p^r}f_j-{x_j}^{p^r}f_i|_K\in H^0\mathscr
I_{Z}(s+p^r)
\]
where $p^r\le m+1$. By the assumption that $f|_K$ and $g$ generate
$I_Z$ in degree $\le s+p^r$ we can
say that:
\[
{x_i}^{p^r}f_j-{x_j}^{p^r}f_i|_K=h_{ij}f|_K+l_{ij}g,
\]
for some $h_{ij},l_{ij}\in H^0\mathscr O_K(p^r)$. So:
\begin{equation}
  \label{eq:4}
  E\cap K=V\left(f|_K,l_{ij}g \mid i,j=0,\dots,4\right).
\end{equation}
So $E\cap K$ contains the complete intersection of 2 curves of degree
$s$ $V(f|_K,g)$, but we have seen that $\deg(E\cap K)\le s^2$. This
implies that $E\cap K$ is the complete intersection $V(f|_K,g)$ and so
$\mathscr I_{E\cap K|R\cap K}\cong \mathscr O_{R\cap K}(-s)$. So by
\eqref{eq:5} we have:
\begin{equation}
  \label{eq:6}
  0\rightarrow \mathscr N|_K\rightarrow {F^r}^{\star}\Omega_K(p^r)\oplus
\mathscr O_K\rightarrow \mathscr O_{R\cap K}\rightarrow 0.
\end{equation}
By the fact that $h^0\mathscr N|_K=0$, that $R\cap K$ is integral and by the commutative diagram:
 \begin{center}
\begin{tikzpicture}
\matrix(m)[matrix of math nodes, row sep=2.6em, column sep=2.8em, text
height=1.5ex, text depth=0.25ex] {& & 0 & &\\ 
& 0  & \mathscr O_K & \mathscr O_K & 0 &\\ 
0 & \mathscr{N}|_K & {F^r}^{\star}\Omega_K(p^r)\oplus \mathscr O_K &
\mathscr O_{R\cap K} & 0\\
& \mathscr N|_K & {F^r}^{\star}\Omega_K(p^r) & &\\
& 0 & 0  & &\\}; 
\path[->,font=\scriptsize]
(m-1-3) edge (m-2-3)
(m-2-2) edge (m-2-3)
(m-2-3) edge (m-2-4)
(m-2-4) edge (m-2-5)
(m-2-2) edge (m-3-2)
(m-2-3) edge (m-3-3)
(m-2-4) edge (m-3-4)
(m-3-1) edge (m-3-2)
(m-3-2) edge (m-3-3)
(m-3-3) edge (m-3-4)
(m-3-4) edge (m-3-5)
(m-3-2) edge (m-4-2) 
(m-3-3) edge (m-4-3) 
(m-4-2) edge (m-4-3) 
(m-4-2) edge (m-5-2) 
(m-4-3) edge (m-5-3); 
\end{tikzpicture}
\end{center}

we get the exact sequence:
\begin{equation}
  \label{eq:7}
  0\rightarrow \mathscr O_K(-s)\rightarrow \mathscr N|_K\rightarrow
  {F^r}^{\star}\Omega_K(p^r) \rightarrow 0.
\end{equation}
By the exact sequence:
\[
0\rightarrow {F^r}^{\star}\Omega_K(p^r)\rightarrow \mathscr O_K^{\oplus
  3}\rightarrow \mathscr O_K(p^r)\rightarrow 0
\]
and by the fact that $p^r\le m+1<\tfrac{s}{2}$ we see that
$\operatorname{Ext}^1({F^r}^{\star}\Omega_K(p^r),\mathscr O_K(-s))=0$
and so $\mathscr N|_K\cong
{F^r}^{\star}\Omega_K(p^r)\oplus \mathscr O_K(-s)$. Since
${F^r}^{\star}\Omega_K(p^r)$ is stable and:
\[
\mu({F^r}^{\star}\Omega_K(p^r))=-\dfrac{p^r}{2}>\mu(\mathscr O_K(-s))=-s,
\]
we see that the maximal destabilizing subsheaf of $\mathscr N|_K$ is
${F^r}^{\star}\Omega_K(p^r)$. So, since $\mathscr N$ is unstable of
rank 3, by \cite[Theorem 3.1]{MARU} the maximal destabilizing  subsheaf of $\mathscr N$ must
be a reflexive sheaf $\mathscr F$ of rank 2 such that:
\begin{equation}
  \label{eq:10}
  \mathscr F|_K\cong  {F^r}^{\star}\Omega_K(p^r).
\end{equation}
So, being $\mathscr F$ the maximal destabilizing sheaf of $\mathscr N$, we have the following commutative diagram:
 \begin{center}
\begin{tikzpicture}
\matrix(m)[matrix of math nodes, row sep=2.6em, column sep=2.8em, text
height=1.5ex, text depth=0.25ex] {& 0 & 0 & &\\ 
0 & \mathscr{F} & \mathscr{F} & 0 &\\ 
0 & \mathscr{N} & {F^r}^{\star}\Omega_H(p^r) & \mathscr I_{E|R}(s) & 0\\
& \mathscr I_T(-s) & \mathscr Q & \mathscr I_{E|R}(s) &\\
& 0 & 0  & 0 &\\}; 
\path[->,font=\scriptsize]
(m-1-2) edge (m-2-2)
(m-1-3) edge (m-2-3)
(m-2-1) edge (m-2-2)
(m-2-2) edge (m-2-3)
(m-2-3) edge (m-2-4)
(m-2-2) edge (m-3-2)
(m-2-3) edge (m-3-3)
(m-2-4) edge (m-3-4)
(m-3-1) edge (m-3-2)
(m-3-2) edge (m-3-3)
(m-3-3) edge (m-3-4)
(m-3-4) edge (m-3-5)
(m-3-2) edge (m-4-2) 
(m-3-3) edge (m-4-3) 
(m-3-4) edge (m-4-4) 
(m-4-2) edge (m-5-2) 
(m-4-3) edge (m-5-3) 
(m-4-4) edge (m-5-4); 
\end{tikzpicture}
\end{center}

where $\mathscr I_T$ is the ideal sheaf in $H$ of a zero-dimensional
scheme $T$ and $\mathscr Q$ is a rank 1 sheaf such that $c_1(\mathscr
Q)=0$. Since $Q|_K\cong \mathscr O_K$,  $\mathscr Q$ must be torsion
free and so  $\mathscr Q=\mathscr I_W$ for some zero-dimensional
scheme $W$. So we get:
\[
    0\rightarrow \mathscr I_T(-s)\rightarrow \mathscr I_W\rightarrow
\mathscr I_{E|R}(s)\rightarrow 0,
\]
by which we get that $W\ne \emptyset$, because $h^0\mathscr
I_Y(s)=1$. Moreover:
\begin{equation}
  \label{eq:15}
  h^1\mathscr I_E(n)=h^1\mathscr I_{E|R}(n)=\deg W-\deg T
\end{equation}
for any $n<s$ and:
\begin{equation}
  \label{eq:16}
  h^1\mathscr I_E(s)=h^1\mathscr I_{E|R}(s)=\deg W-\deg T-1,
\end{equation}
because $h^0\mathscr I_{E|R}(s)=0$.

Let $F\subset E$ be the equidimensional component of dimension $1$. Then there
exists a sheaf $\mathscr K$ of finite length determining the following
exact sequence:
\[
0\rightarrow \mathscr I_E\rightarrow \mathscr I_F\rightarrow \mathscr K\rightarrow 0.
\]
Then we see that $h^1\mathscr I_E(n)=h^0\mathscr K$ for $n\ll 0$, so
that by \eqref{eq:15} we see that $h^0\mathscr K=\deg W-\deg T$. Moreover:
\[
h^0\mathscr I_E(s)-h^0\mathscr I_F(s)+h^0\mathscr K-h^1\mathscr
I_{E}(s)+h^1\mathscr I_F(s)=0
\]
and so, since $Y\subset F\subset E$, $h^0\mathscr I_E(s)=h^0\mathscr I_F(s)=1$ and by \eqref{eq:16} we get:
\[
h^1\mathscr I_{F}(s)=h^1\mathscr I_E(s)-h^0\mathscr K=-1.
\]
This is impossible and so we get a contradiction. 
\end{proof}

\begin{corollary}
  \label{C:1}
  Let $h^0\mathscr I_Y(s)\ne 0$ and let $p<s$. If $\deg X>s^2$,
  then $h^0\mathscr I_X(s)\ne 0$.
\end{corollary}

In the following theorem we see that for $p\ge s$ the bound for $d$ is
independent on the the order of the sporadic zero $\alpha$ and
coincides with the bound of the characteristic zero case (see \cite{MeR} and \cite{Me1}).

\begin{theorem}
  \label{T:2}
  Let $h^0\mathscr I_Y(s)\ne 0$, $h^0\mathscr I_X(s)=0$ and let $p\ge
  s$. Then $\deg X\le s^2-s+2$.
\end{theorem}
\begin{proof}
  The proof works as in Theorem \ref{T}. We just need to remark that
  in the case $p\ge s$ it must be $r=0$, which means $p^r=1$. Indeed,
  proceeding as in Lemma \ref{L:1} we see that we get an exact
  sequence:
\[
0\rightarrow \mathscr I_{X}(s-p^r)\rightarrow \mathscr
I_X(s)\rightarrow \mathscr I_{X\cap H_r|H_r}(s)\rightarrow 0,
\]
where $\mathscr I_{X\cap H_r|H_r}\subset \mathscr O_{H_r}$ is the
ideal sheaf of $X\cap H_r$. Since $h^0\mathscr I_{X\cap H_r|H_r}(s)\ne
0$ and $h^0\mathscr I_{X}(s)=0$, it must be $h^1\mathscr I_X(s-p^r)\ne
0$. By the fact that $X$ is integral we see that it must be $p^r<s$
and so $r=0$ and $p^r=1$.
\end{proof}

Now we show that the bounds given in Theorem \ref{T} and Theorem
\ref{T:2} are sharp.

\begin{example}  \label{ex}
  Let $r$, $p$, $s\in \mathbb N$ such that $s\ge 2p^r$. Let us
  consider $\mathscr
  E=\mathscr O_{\mathbb P^4}(p^r-2s)\oplus \mathscr O_{\mathbb
    P^4}(-p^r-s) ^{\oplus 2}$ and $\mathscr
  F={F^r}^{\star}\Omega_{\mathbb P^4}(p^r-s)$. Then, since
  $E^{\vee}\otimes F$ is generated by global sections, by \cite{K} the
  dependency locus of a general homomorphism $\varphi\in
  \operatorname{Hom}(\mathscr E,\mathscr F)$ is a smooth
  surface $X\subset \mathbb P^4$ and it is determined by the sequence:
   \begin{equation}
    \label{eq:14}
    0\rightarrow \mathscr O_{\mathbb P^4}(p^r-2s)\oplus \mathscr O_{\mathbb
    P^4}^{\oplus 2}(-p^r-s)\rightarrow {F^r}^{\star}\Omega_{\mathbb
    P^4}(p^r-s)\rightarrow \mathscr I_X\rightarrow 0.
  \end{equation}
Together with:
\begin{equation}
  \label{eq:13}
  0\rightarrow {F^r}^{\star}\Omega_{\mathbb
    P^4}(p^r)\rightarrow \mathscr O_{\mathbb P^4}\rightarrow \mathscr
  O_{\mathbb P^4}(p^r)\rightarrow 0
\end{equation}
this implies that $h^1\mathscr I_X=0$, so that $h^0\mathscr O_X=1$ and
$X$ is connected and,
being smooth, $X$ is integral. Moreover, $h^0\mathscr I_X(s)=0$ and by
a computation with Chern classes we see that $\deg X=s^2-p^rs+2p^{2r}$. 

Let $H\subset \mathbb P^4$ be a general
hyperplane and let $H_r\subset \mathbb P^4$ be the nonreduced hypersurface of degree
$p^r$ such that $H_r|_{\text{red}}=H$. Then, $(F^r)^{-1}(H)=H_r$. This
shows that we have a commutative diagram:
\begin{center}
\begin{tikzpicture}
\matrix(m)[matrix of math nodes, row sep=2.6em, column sep=2.8em, text
height=1.5ex, text depth=0.25ex] {H_r&H\\ \mathbb P^4&\mathbb P^4\\}; 
\path[->] (m-1-1) edge node[auto] {$\scriptstyle \pi$} (m-1-2);
\path[right hook->] (m-1-1) edge node[auto, left] {$\scriptstyle i$} (m-2-1) ;
\path[right hook->] (m-1-2) edge node[auto] {$\scriptstyle j$} (m-2-2); 
\path[->] (m-2-1) edge node[below] {$\scriptstyle{F^r}$} (m-2-2); 
\end{tikzpicture}
\end{center}
So we have:
\[
i^{\star}({F^r}^{\star}\Omega_{\mathbb
  P^4}(p^r))=i^{\star}({F^r}^{\star}(\Omega_{\mathbb
  P^4}(1)))={\pi^{\star}}(j^{\star}(\Omega_{\mathbb P^4}(1)))\cong
{\pi^{\star}}(\Omega_H(1))\oplus \mathscr O_{H_r}.
\]
This implies that $h^0({F^r}^{\star}\Omega_{\mathbb
  P^4}(p^r)|_{H_r})\ge 1$. In particular, by \eqref{eq:14} we see that
$h^0\mathscr I_{X\cap H_r|H_r}(s)\ne 0$, so that $h^0\mathscr
I_Y(s)\ne 0$. Moreover, by \eqref{eq:14} and by \eqref{eq:13}
we see that $h^1\mathscr
I_X(s-p^r-1)=0$. This shows that $X$ has a sporadic zero of degree $s$
and order $m=p^r-1$. So:
\begin{enumerate}
\item if $r=0$ and $s\ge 2$, then $p^r=1$, $m=0$ and $\deg X=s^2-s+2$;
\item if $s=2p^r+1$, then $s=2m+3$ and $\deg X=s^2-\tfrac{s-1}{2}=s^2-s+p^r+1$;
\item if $s=2p^r$, then $s=2m+2$ and $\deg X=s^2$.
\end{enumerate}
This shows that the bounds in Theorem
  \ref{T} and Theorem \ref{T:2} are sharp.
\end{example}

\end{document}